\documentclass{article} 
\usepackage{hyperref}
\usepackage[accepted]{icml2023}
\usepackage{times}
\usepackage{tabularx}
\usepackage{graphicx}
\usepackage{amsmath}
\usepackage{amssymb}
\usepackage{booktabs}
\usepackage[utf8]{inputenc}
\usepackage{algorithm}
\usepackage{caption}
\usepackage{subcaption}
\usepackage{graphicx}
\usepackage{array}
\usepackage{enumitem}
\setlist{nosep}

\usepackage{amsmath,amsfonts,bm}

\def\eqref#1{equation~\ref{#1}}

\def\1{\bm{1}}

\DeclareMathAlphabet{\mathsfit}{\encodingdefault}{\sfdefault}{m}{sl}
\SetMathAlphabet{\mathsfit}{bold}{\encodingdefault}{\sfdefault}{bx}{n}

\newcommand{\R}{\mathbb{R}}

\usepackage{url}
\usepackage{multirow}

\newcommand{\revised}[2]{{#1}}

\newcommand{\bb}{{\boldsymbol b}}

\newcommand{\pb}{{\boldsymbol p}}

\newcommand{\rbd}{{\boldsymbol r}}

\newcommand{\wb}{{\boldsymbol w}}
\newcommand{\xb}{{\boldsymbol x}}

\newcommand{\zb}{{\boldsymbol z}}

\newcommand{\pd}[2]{ \dfrac{\partial #1}{\partial #2}}

\newcommand{\pdtwo}[2]{ \dfrac{\partial^2 #1}{\partial #2^2}}

\newcommand{\grad}{{\nabla}}

\begin{document}

\twocolumn[
\icmltitle{Learning Preconditioners for Conjugate Gradient PDE Solvers}



\icmlsetsymbol{equal}{*}

\begin{icmlauthorlist}
\icmlauthor{Yichen Li}{mitcsail}
\icmlauthor{Peter Yichen Chen}{mitcsail}
\icmlauthor{Tao Du}{thu,sh}
\icmlauthor{Wojciech Matusik}{mitcsail}
\end{icmlauthorlist}

\icmlaffiliation{mitcsail}{MIT CSAIL}
\icmlaffiliation{thu}{Tsinghua University}
\icmlaffiliation{sh}{Shanghai Qi Zhi Institute}

\icmlcorrespondingauthor{Yichen Li}{yichenl@csail.mit.edu}

\icmlkeywords{Machine Learning, ICML}

\vskip 0.3in
]

\printAffiliationsAndNotice 



%


\begin{abstract}
Efficient numerical solvers for partial differential equations empower science and engineering.
One commonly employed numerical solver is the preconditioned conjugate gradient (PCG) algorithm, whose performance is largely affected by the preconditioner quality. However, 
designing high-performing preconditioner with traditional numerical methods is highly non-trivial, often requiring problem-specific knowledge and meticulous matrix operations. 
We present a new method that leverages learning-based approach to obtain an approximate matrix factorization to the system matrix to be used as a preconditioner in the context of PCG solvers. 
Our high-level intuition comes from the shared property between preconditioners and network-based PDE solvers that excels at obtaining approximate solutions at a low computational cost. Such observation motivates us to represent preconditioners as graph neural networks (GNNs). In addition, we propose a new loss function that rewrites traditional preconditioner metrics to incorporate inductive bias from PDE data distributions, enabling effective training of high-performing preconditioners. 
We conduct extensive experiments to demonstrate the efficacy and generalizability of our proposed approach on solving various 2D and 3D linear second-order PDEs.\footnote{\url{https://sites.google.com/view/neuralPCG}}{}


\end{abstract}

\section{Introduction}

 The conjugate gradient (CG) algorithm is an efficient numerical method for solving large sparse linear systems. CG iteratively reduces the residual error to solve the linear systems to a specified accuracy level and does not require the expensive computation of a full matrix inverse. CG is commonly applied to solve the underlying large sparse linear systems originating from discretized partial differential equations (PDEs), as PDEs generally lack analytical, closed-form solutions. 
 Therefore, developing a high-quality linear solver is instrumental in solving numerical PDEs efficiently and effectively. 

A naive CG implementation suffers from a slow convergence rate for ill-conditioned matrices. Therefore, it is typically equipped with a \emph{preconditioner} that modulates the system matrix's condition number. Preconditioners are mathematically-grounded and problem-dependent transformation matrices that can be applied to the original linear system $\mathbf{A}\mathbf{x}=\mathbf{b}$, where $A$ is a sparse matrix, $\mathbf{A} \in \mathbb{R}^{n\times n}$, and $\mathbf{x}, \mathbf{b} \in \mathbb{R}^{n}$ are the unknown solution vector and right-hand-side vector respectively. Preconditioner $\mathbf{P} \in \mathbb{R}^{n\times n}$ transforms the original into an easier one when applied to both sides of the original linear system, $\mathbf{P}^{-1}\mathbf{Ax}=\mathbf{P}^{-1}\mathbf{b}$. Traditional numerical preconditioners use the approximation of $\mathbf{A}$ as preconditioners, and the quality of preconditioning significantly affects the CG solvers' convergence rate.

Designing a high-performance preconditioner is challenging because an ideal preconditioner is inexpensive to solve and greatly reduces the condition number of the system matrix. These two desired properties that often conflict. 
Decades of research effort in applied math have been devoted to addressing this longstanding issue through proposing new sparsity patterns~\cite{schafer2021sparse, chen_2021_icsiggraph, davis_sparse_cholesky}, matrix reordering techniques~\cite{Liu_book, chen2021bals, brandhorst2011fast}, and multi-level approaches~\cite{SAAD200199_multilevel, buranay2019approximate, liew2007computational}. 
Recent efforts have been made towards leveraging machine learning and optimization to discover preconditioners~\cite{gotz_arxiv_cnnprecond, sappl_arxiv, akmann_arxiv, azulay_arxiv}, but their methods are typically tailored to specific PDE problems.
To our best knowledge, a learning-based solution applicable to general PDE problems is unavailable, and our work aims to fill this gap.

In this work, we propose a learning-based method for preconditioning sparse symmetric positive definite (SPD) matrices. Recent machine learning works on simulating PDE-governed physical systems~\cite{pfaff2020learning, sanchez2020learning} demonstrate the efficacy of neural networks in obtaining cheap approximations of PDE solutions with modest accuracy. 
The high-level intuition is that this nature aligns well with the desired role of a preconditioner.
With this intuition, we propose to learn an approximate decomposition of the system matrix itself to be used as a preconditioner in the context of CG solvers. 
We leverage the duality between graph and matrix to propose a generic preconditioning solution using the graph neural network (GNN). 
To facilitate the training of our neural network preconditioner, we propose a novel preconditioner metric by rewriting classical metrics with considerations of data distribution. Motivated by the classical projected conjugate gradient algorithm~\cite{Gould2001OnTS}, our new preconditioner metric reveals an interesting interpretation of preconditioners through the lens of the bias-variance trade-off in statistics and machine learning. Specifically, when the problem domain is not full rank, there exist projections of the system matrix to the solution subspace. Compared with existing works~\cite{TrotMult2001, cali_arxiv, sappl_arxiv, akmann_arxiv, azulay_arxiv}
that focus exclusively on the system matrix $\mathbf{A}$, our proposed method also attends to the data distribution of solution vector $\mathbf{x}$. 

Our approach has several benefits: the learned preconditioner outperforms classical preconditioners because it adapts to data distributions from target PDE applications. Additionally, the duality between graph and matrix grants our learned preconditioners the excellent property of correctness-by-construction and order-invariance. 
Moreover, unlike multi-grid preconditioners that specialize in elliptic PDEs~\cite{TrotMult2001} or other learning-based preconditioners~\cite{cali_arxiv, sappl_arxiv, akmann_arxiv, azulay_arxiv} that focus on a single problem or PDE, our proposed method for preconditioning is generic and applicable to many different problems. We demonstrate this property with elliptic, parabolic, and hyperbolic PDEs.

To showcase the efficacy of our proposed preconditioner, we evaluate it on a set of 2D and 3D representative linear second-order PDEs. We compare its performance with (1) existing learning-based preconditioning methods and (2) classical numerical preconditioners. Our experiments show that our proposed method has a speed advantage over both baselines. It learns a preconditioner tailored to the training data distribution, which other preconditioning methods do not exploit. Finally, we also demonstrate the strong generalizability of our proposed approach with respect to varying physical parameter values and geometrical domains. 

In summary, our work makes the following contributions:
\begin{itemize}[leftmargin=*]
\setlength\itemsep{1mm}
    \item We propose a generic learning-based framework to precondition the conjugate gradient algorithm. Our proposed method guarantees positive definiteness by design and works with irregular geometric domains. 
    \item We propose a novel and efficient loss function that introduces inductive bias to preconditioning large and sparse system matrices. 
    \item We conduct extensive experiments to evaluate the efficacy and generalizability of our proposed approach and demonstrate its advantages over existing methods.
\end{itemize}

\section{Related Work}\label{sec:related}

\paragraph{Numerical preconditioning} Preconditioning is a classical numerical technique in solving linear systems of equations. It typically applies a carefully chosen matrix to transform a linear system into one with a smaller condition number. Below, we briefly review representative works from three technical aspects of a numerical preconditioner: matrix factorization~\cite{golub2013matrix, khare2012sparse}, matrix reordering~\cite{Liu_book,davis1999modifying,schafer2021sparse}, and multiscale approaches~\cite{chen2021multiscale}.

Matrix factorization inspires several widely used numerical preconditioners, but they face the problem of speed and accuracy trade-off. For example, in the extreme cases, the Jacobi preconditioner is a direct inverse of the diagonal elements of the original matrix $\mathbf{A}$. it is fast to derive but is very limited in reducing the condition number of the original matrix $\mathbf{A}$. The famous incomplete Cholesky (IC) preconditioner~\cite{nocedal1999numerical} originates from the Cholesky decomposition of SPD matrices. Such preconditioners face the trade-off between speed and accuracy: a complete factorization essentially solves the SPD matrix but is very expensive, whereas an incomplete factorization is cheaper to compute but has a limited impact on improving the condition number. The more advanced factorization approach considers fill-in to capture the additional non-zero entries~\cite{johnson2012numerical}, but it is more computationally expensive to construct, as it requires multiple rounds of factorization.  \citet{schafer2021sparse} address this trade-off using a Kullback-Leibler minimization approach to speed up factorization-based methods. Our proposed method tackles the same challenge using a neural network and leverages the data distribution of the linear systems.

Matrix reordering techniques \cite{Liu_book, chen2021bals, brandhorst2011fast} aim to reduce the matrix bandwidth by reshaping the sparse matrices with large bandwidths to block diagonal form. They are commonly used with other numerical preconditioners, such as factorization-based ones, to reduce fill-in and improve the parallelizability when deriving the preconditioners. Fortunately, our proposed method uses a graph-neural-network (GNN) representation and inherits its excellent property of order invariance and parallelizability by design.

Finally, multi-level approaches help to improve the scalability of a standard numerical preconditioner. A representative multi-level technique is the multigrid method~\cite{TrotMult2001, SAAD200199_multilevel,  chen_2021_icsiggraph}, which uses smoothing to communicate between coarser and finer discretizations. It is a power tool for preconditioning elliptic PDEs but struggles with hyperbolic and parabolic PDEs~\cite{TrotMult2001}. On the contrary, our method follows an orthogonal direction by studying the numerical approximation of a given system and discretization and is not tailored to specific PDE instances.  Research has shown that it is possible to combine the two directions~\cite{SAAD200199_multilevel, chen_2021_icsiggraph}.

\paragraph{Learning-based preconditioning} More recent works on preconditioner design borrows inspirations from machine learning techniques~\cite{belbute2020combining,um2020solver,li2020fourier,li2020multipole,raissi2019physics,karniadakis2021physics}. Similar to our work, several recent papers~\cite{azulay_arxiv, sappl_arxiv, akmann_arxiv, cali_arxiv, koolstra_arxiv} also model preconditioners with neural networks, but their convolutional-neural-network (CNN) architectures are strongly correlated with a grid discretization of a rectangular domain. However, we are different from these works in three folds. We leverage the GNN architecture that is mesh-friendly and applicable to irregular geometrical boundaries. In addition, unlike previous works that uses hard-coded threshold to satisfy the constraint for a specific problem, our proposed method leverage the duality between graph and matrix and works on different representative second-order linear PDEs. Finally, our work differs from these papers in formulating a novel preconditioner metric as the training loss function, which incorporates data distributions often overlooked before into preconditioners.

\section{Preliminaries}\label{sec:pre}
\paragraph{Linear second-order PDEs} We consider linear second-order PDEs in the following format:
\begin{align}
\frac{1}{2}\nabla \cdot \mathbf{K}\nabla f(\mathbf{p}) + \mathbf{a} \cdot \nabla f(\mathbf{p}) = c(\mathbf{p}), \quad \forall \mathbf{p}\in\Omega.
\end{align}
Here, $\Omega\subset\R^d$ ($d=2$ or $3$) is the problem domain, $f:\Omega\rightarrow \mathbb{R}$ is the function to be solved, $\mathbf{K}\in\mathbb{R}^{d \times d}$ and $\mathbf{a}\in\mathbb{R}^d$ are constants, and $c:\Omega\rightarrow\mathbb{R}$ is a given source function. We assume $\mathbf{K}$ to be symmetric, whose eigenvalues classify these PDEs into \emph{elliptic} (e.g., the Poisson or Laplace equation), \emph{hyperbolic} (e.g., the wave equation), and \emph{parabolic} (e.g., the heat equation) equations.

\paragraph{Boundary conditions} We equip the PDEs with \emph{Neumann} and \emph{Dirichlet} boundary conditions:
\begin{align}
\frac{\partial f(\mathbf{p})}{\partial \mathbf{n}} =& N(\mathbf{p}), \quad \forall \mathbf{p}\in\partial\Omega_N, \\
f(\mathbf{p}) = & D(\mathbf{p}), \quad \forall \mathbf{p}\in\partial \Omega_D.
\end{align}
Here, $\partial \Omega_N$ and $\partial \Omega_D$ form a partition of the domain boundary $\partial \Omega$, and $N:\partial\Omega_N\rightarrow \mathbb{R}$ and $D: \partial\Omega_D\rightarrow\mathbb{R}$ are two user-specified functions. The notation $\frac{\partial}{\partial \mathbf{n}}$ represents the directional derivative along the normal ($\mathbf{n}$) direction.

\paragraph{Discretization} PDEs are continuous problems that must be discretized before applying a numerical solver. 
We adopt the standard Galerkin method from the finite element theory \citep{johnson2012numerical}, resulting in a linear system
\begin{align}
\mathbf{A}\mathbf{x} = \mathbf{b},
\end{align}
where $\mathbf{A}\in\mathbb{R}^{n \times n}$, with $n$ the number of degrees of freedom (DoFs) after discretization, is the \emph{stiffness matrix} of the PDE system, which is sparse and SPD. The vector $\mathbf{b}\in\mathbb{R}^n$ typically contains information from the source term and the (discretized) boundary conditions. The goal is to solve for $\mathbf{x}\in\mathbb{R}^n$, which approximates the field of interest $f$ at discretized locations of DoFs in $\Omega$.

\paragraph{PCG algorithm} The PCG algorithm takes a system matrix $\mathbf{A}$ and a right-hand side vector $\mathbf{b}$ to solve for $\mathbf{x}$. It starts with an initial guess $\mathbf{x}_0$ and iteratively updates it by moving towards conjugate directions for suppressing the residual $\mathbf{r}=\mathbf{b}-\mathbf{Ax}$. The preconditioner $\mathbf{P}$ transforms the original problem such that the gradient direction is from the preconditioned system residual $\mathbf{z} = \mathbf{P}^{-1}\mathbf{r}$, as shown in Alg.~\ref{alg:pcg} in Appendix.  

The condition number indicates the extent to which the preconditioner can reduce the number of CG iterations, 
$$ \kappa(\mathbf{A}) = \frac{\lambda_{\max} (\mathbf{A})}{ \lambda_{\min} (\mathbf{A})},$$ 
where $\lambda_{\max}(\mathbf{A})$ and $\lambda_{\min}(\mathbf{A})$ denote the largest and smallest eigenvalue of $\mathbf{A}$, respectively. A good preconditioner $\mathbf{P}$ can transform the original system $\mathbf{A}$ into an easy-to-solve matrix with clustered eigenvalues. Alg.~\ref{alg:pcg} reveals that a high-performing $\mathbf{P}$ needs to be an easily invertible SPD matrix that is similar to $\mathbf{A}$.

\section{Method}\label{sec:method}

\subsection{Problem Setup}\label{sec:method:setup}

\label{sec:network-precond}

Our task is to learn a mapping $\mathbf{P} = f_\theta (\mathbf{A}, \mathbf{b})$, where $\mathbf{A}$ is a sparse SPD matrix derived from a given PDE problem, $\mathbf{b}$ is the right-hand-side vector, $\theta$ are the learned parameters, and $\mathbf{P}$ is the resulting preconditioning matrix.
We leverage the duality between graph and matrix to ensure that our learned $\mathbf{P}$ is valid (SPD and easily-invertible).

\subsection{Learning Preconditioners}

The main consideration for choosing a graph neural network over a convolutional neural network for predicting preconditioner~\cite{sappl_arxiv, azulay_arxiv} comes from the duality between graphs and square matrices. Graphs are composed of a set of nodes and edges $(\{v_{i}\}, \{e_{i,j}\})$ . Each node can correspond to a row or column in the matrix, and each edge corresponds to an entry in the matrix. In a sparse matrix, the corresponding graph edges only run between the nonzero entries in the matrix. 
The nice duality also directly transfers to PDE problems on $\Omega$ discretized as a triangle mesh in 2D or a tetrahedron mesh in 3D. Graph nodes and edges directly correspond to mesh vertices and edges, respectively.

We store the input $\mathbf{A}$ as a one-dimensional edge feature: If $\mathbf{A}_{ij}$ is a nonzero entry in $\mathbf{A}$, we add it as an edge feature connecting node $i$ to $j$. Similarly, we store vector $\mathbf{x}$ as a one-dimensional node feature on the graph. A matrix-vector product such as $\mathbf{A}\mathbf{x}$ can be viewed as a round of message passing on the graph. Each node sends a message to its neighboring nodes, which in turn pass the message along to their own neighbors. The message at each node is updated based on the messages it receives from its neighbors, which can be represented by a weighted sum of the messages. This operation gives us a new node value corresponding to the resulting vector from the multiplication.

\paragraph{Architecture} We use a variant of the encoder-processor-decoder architecture from previous literature~\cite{pfaff2020learning, luz2020amg, sanchez2020learning}. There are three main components to this architecture. The encoder uses an MLP, which takes as input the graph nodes and edges and outputs a 16D feature. The feature representations are updated through a series of MLP message-passing layers in the processing stage, where nodes and edge features are updated through aggregating features in the local neighborhood. We use five message-passing layers in the processors. Finally, the decoder takes the updated feature on each edge to predict a real-number value on each edge which forms a matrix $\mathbf{M}$ to be used to construct our predicted preconditioner $\mathbf{P}$.

Unfortunately, directly assembling the predicted edge feature into a matrix often fails to serve as a valid preconditioner because there is no guarantee of its symmetry or positive definiteness. Therefore, we first construct a triangular matrix by averaging a pair of the bidirectional edges running between the two graph nodes and store the value on the corresponding lower-triangular indices $\mathbf{L}_{i, j | i \geq j }$. 
We use $\mathbf{LL}^\top$ as our preconditioner, and this construction ensures its symmetry and positive definiteness. \revised{Please see Appendix Section~\ref{supp:sec:ours} for more details on the GNN architecture.}

\subsection{Loss function}

Existing works on learning preconditioners~\cite{sappl_arxiv, cali_arxiv, azulay_arxiv, akmann_arxiv} leverage the loss function that minimizes the condition number $\kappa$ of the system matrix $\mathbf{A}$. Condition number can be a natural choice when designing loss functions to discover new preconditioners since the strength of preconditioners in reducing CG convergence iterations is greatly reflected by condition number $\kappa$. However, using condition number as a loss function has two main drawbacks: 1) condition number is expensive and slow to compute because every data instance of $\mathbf{A}$ requires a full eigen decomposition, which is a $O(n^3)$ operation, making it very restrictive in training for large system of equations. 2) Previous literature~\cite{wang2019backpropagation} has reflected that backpropagation through eigen decomposition tends to be numerically unstable. Therefore, we relax the objective to the squared Frobenius norm. 

Since our formulation uses matrix decomposition for approximating the original system matrix $\mathbf{A}$, designing many classic preconditioners can be cast as a problem of minimizing their discrepancy to the given linear system over a set of easy-to-compute matrices:
\begin{align}
\min_{\mathbf{P}\in\mathcal{P}} L(\mathbf{P}, \mathbf{A}),
\end{align}
where $\mathbf{A}$ is the system matrix defined above and the system that we want to precondition upon, $\mathcal{P}$ is the feasible set of preconditioners, and $L(\cdot, \cdot)$ is a loss function defined on the difference between the two input matrices. 

The design of classic preconditioners, e.g., incomplete Cholesky or symmetric successive over-relaxation (SSOR) \citep{golub2013matrix}, is defined on the left-hand-side matrix $\mathbf{A}$ only. 
Following this classical approach, it is now tempting to consider the following loss function definition for our neural-network preconditioner:
\begin{align}
\label{eqn:loss_A}
\min_{\mathbf{\theta}} \sum_{(\mathbf{A}_i, \mathbf{x}_i, \mathbf{b}_i)} \| \mathbf{L}_\mathbf{\theta}(\mathbf{A}_i, \mathbf{b}_i)\mathbf{L}^\top_\mathbf{\theta}(\mathbf{A}_i, \mathbf{b}_i) - \mathbf{A}_i \|^2_F,
\end{align}
where $\mathbf{\theta}$ is the network parameters to be optimized, $\mathbf{L}(\mathbf{\theta})$ is the lower-triangular matrix from the network's output, and $\|\cdot\|_F^2$ represents the squared Frobenius norm. The index $i$ loops over training data tuples $(\mathbf{A}_i, \mathbf{x}_i, \mathbf{b}_i)$. This definition closely resembles the goal of the famous incomplete Cholesky preconditioner, especially since $\mathbf{L}$ shares the same sparsity pattern as the lower triangular part of $\mathbf{A}$.

We argue that this design decision unnecessarily limits the full power of preconditioners because they overlook the right-hand-side vector $\mathbf{b}$ and its distribution among actual PDE problem instances. A closer look at the loss function can reveal the potential inefficiencies in its design:
\begin{align}
L := & \sum_i \|\mathbf{L}_\theta \mathbf{L}_\theta^\top - \mathbf{A}_i \|_F^2 \\
= & \sum_i \|(\mathbf{L}_\theta \mathbf{L}_\theta ^\top - \mathbf{A}_i)\mathbf{I} \|_F^2 \\
= & \sum_i \sum_j \|\mathbf{L}_\theta \mathbf{L}_\theta ^\top\mathbf{e}_j  - \mathbf{A}_i\mathbf{e}_j\|_F^2
\end{align}
where $\mathbf{e}_j$ stands for the one-hot vector with one at the $j$-th entry and zero elsewhere. This derivation shows that this loss encourages a well-rounded preconditioner with uniformly small errors in all $\mathbf{e}_j$ directions, regardless of the actual data distribution in the training data $(\mathbf{A}_i, \mathbf{x}_i, \mathbf{b}_i)$. 

In contrast to these classic preconditioners, we propose to learn a neural network preconditioner from both left-hand-side matrices and right-hand-side vectors in the training data. 
Therefore, we consider a new loss function instead:
\begin{align}
\label{eqn:loss_Ab}
L := & \sum_i \|\mathbf{L}_\theta \mathbf{L}_\theta^\top\mathbf{x}_i - \mathbf{A}_i\mathbf{x}_i \|_2^2 \\
= & \sum_i \|\mathbf{L}_\theta \mathbf{L}_\theta^\top\mathbf{x}_i - \mathbf{b}_i \|_2^2.
\end{align}
Comparing these two losses, we can see that the new loss replaces $\mathbf{e}_j$ with $\mathbf{x}_i$ from the training data. Therefore, the new loss encourages the preconditioner to ensemble $\mathbf{A}$ not uniformly in all directions but towards frequently seen directions in the training set. Essentially, this new loss trades generalization of the preconditioner with better performance for more frequent data.

\subsection{Remark}
To summarize, we overcome the traditional trade-off between speed and efficacy by carefully limiting the scope of our preconditioners by leveraging the speed and approximation ability of neural networks. Our approach leverages the duality between graph and matrix to guarantee correctness-by-construction (ensure that the learned preconditioner is an SPD matrix). This construction also allows for application to various PDE problems. Finally, we also exploit the data distribution in our loss function design for more efficient training and learning of a more effective preconditioner.

\revised{There are two directions for speeding up PCG solvers. One is by using an efficient and effective precondition method, as shown in our previous discussion; another direction comes from using a better starting guess $\mathbf{x}_0$. Our proposed method can achieve both simultaneously. In addition to a preconditioner, our method can also function as a surrogate model to offer further speed up without additional computation time cost. This can be achieved by predicting an initial guess $\hat{\mathbf{x}}_0$ for the conjugate gradient algorithm. We directly regress decoded graph node values to the solution $\mathbf{x}_i$ of the linear systems $\mathbf{A}_i \mathbf{x}_i = \mathbf{b}_i$. The predicted $\hat{\mathbf{x}}_i$ can be used as the initial starting point of the CG algorithm. In the experiments section~\ref{sec:exp}, we show results without using the predicted $\hat{\mathbf{x}}_0$ to focus solely on the effect of preconditioners, and we conduct additional experiments that show the additional speed up in Section.~\ref{supp:x-compare}. in Appendix.}

\begin{figure*}
    \centering
    \includegraphics[width=\textwidth]{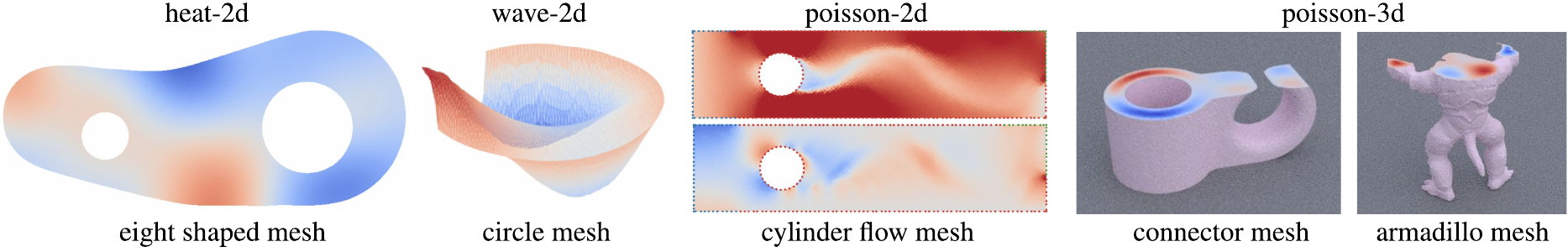}
    \vspace{-4mm}
    \caption{{Environment Overview. Left to right: {heat-2d}, {wave-2d}, {poisson-2d}, and {poisson-3d}.}}
    \label{fig:env}
\end{figure*}

\begin{table*}[htb]
\resizebox{\textwidth}{!}{%
        \scriptsize
\centering
\begin{tabular}{lcccccccc}

\toprule
Task & Method    & Precompute time 	$\downarrow$	& time (iter.) $\downarrow$ 		& time (iter.) 	$\downarrow$	& time (iter.) 	$\downarrow$	& time (iter.) $\downarrow$		& time (iter.) 	$\downarrow$	& time (iter.)$\downarrow$ \\ &

& (s) 		&  until 1e-2 		&  until 1e-4 		&  until 1e-6 		&  until 1e-8 		&  until 1e-10 		&  until 1e-12 \\

\midrule
\multirow{4}{*}{heat-2d} &

 Jacobi & \textbf{0.0001} &  0.657  (32) &  2.188  (132) &  3.263  (202) &  4.105  (257) &  5.269  (333) &  6.255  (398) \\

&Gauss-Seidel & 0.0071 &  0.645  (27) &  1.656  (98) &  2.339  (146) &  2.995  (193) &  3.771  (247) &  4.402  (292) \\

& IC& 1.5453 &  1.954  ({12}) &  2.612  ({54}) &  3.061  ({83}) &  3.409  ({105}) &  3.832  ({133}) &  4.271  ({161}) \\

&IC(2) & 2.3591 &  2.624  (11) &  3.094  (40) &  3.399  (60) &  3.651  (76) &  3.965  (96) &  4.260  (115) \\
&Ours & 0.0251 & \textbf{0.490} (17) & \textbf{1.284} (71) & \textbf{1.856} (110) & \textbf{2.300} (140) & \textbf{2.831} (177) & \textbf{3.377} (214) \\

\midrule
\multirow{4}{*}{wave-2d}
& Jacobi & \textbf{0.0001} &  0.141  (0) &  0.141  (0) &  0.141  (0) &  0.176  (6) &  0.295  (26) &  0.417  (46) \\

&Gauss-Seidel & 0.0079 &  0.089  (0) &  0.089  (0) &  0.09  (0) &  0.1  (2) &  0.16  (11) &  0.232  (22) \\

& IC& 0.7679 &  0.885  ({0}) &  0.885  ({0}) &  0.885  ({0}) &  0.904  ({3}) &  0.953  ({11}) &  1.007  ({20}) \\

&IC(2) & 1.1831 &  1.226  (0) &  1.226  (0) &  1.226  (0) &  1.266  (5) &  1.326  (12) &  1.385  (18) \\

&Ours & 0.0147 & \textbf{0.081} (0) & \textbf{0.081} (0) & \textbf{0.081} (0) & \textbf{0.100} (3) & \textbf{0.156} (12) & \textbf{0.211} (21) \\

\midrule
\multirow{4}{*}{possion-2d} &

 Jacobi & \textbf{0.0001} & 0.980 (275) &1.231 (348) &1.572 (448) &1.822 (522) &2.119 (611) &2.405 (697)\\

&Gauss-Seidel & 0.0071 &  0.699  (194) &  0.964  (273) &  1.26  (361) &  1.518  (438) &  1.807  (525) &  2.099  (613) \\

& IC& 0.7093 &  1.188  ({135}) &  1.309  ({171}) &  1.468  ({219}) &  1.559  ({246}) &  1.774  ({311}) &  1.900  ({349}) \\

&IC(2) & 1.205 &  1.308  (60) &  1.439  (74) &  1.543  (100) &  1.604  (115) &  1.664  (131) &  1.747  (151) \\

&Ours & 0.0145 & \textbf{0.639} (175) &\textbf{0.818} (227) &\textbf{1.017} (286) &\textbf{1.118} (316) &\textbf{1.312} (374) &\textbf{1.510} (432)\\

\midrule
\multirow{4}{*}{possion-3d}

& Jacobi & \textbf{0.0002} & \textbf{1.526} (0) & \textbf{2.693} (7) & 5.496 (25) & 9.552 (50) & 13.636 (76) & 17.080 (97) \\

&Gauss-Seidel & 0.3381 &  5.824  (0) &  6.775  (6) &  9.074  (19) &  12.305  (38) &  15.454  (56) &  18.333  (72) \\
&IC & 9.6878 &  10.668  (1) &  11.353  (6) &  12.592  (15) &  13.826  (23) &  14.954  (31) &  15.812  (37) \\
&IC(2) & 17.138 &  18.083  (1) &  18.661  (5) &  19.667  (11) &  20.599  (17) &  21.704  (24) &  22.560  (30) \\
&Ours & 0.4137 & {3.010} (0) & {3.220} (2) & \textbf{4.815} (13) & \textbf{6.908} (28) & \textbf{8.749} (41) & \textbf{10.406} (53) \\

\bottomrule

\end{tabular}
}
\vspace{-2mm}
\caption{Comparison between preconditioners with PCG. We report precompute time, total time ( ICl. precompute time) for each precision level, and the PCG iterations (in parenthesis). The best value is in bold. $\downarrow$: the lower the better.}
\label{tab:pcg}
\vspace{-2mm}
\end{table*}

\section{Experiments}\label{sec:exp}

Our experiments aim to answer the following questions:

1. How does the proposed method compare with classical and learning-based preconditioners in speed and accuracy?

2. Is our data-dependent loss function effective?

3. Does our approach generalize well to unseen inputs?


We introduce the experiment setup in Sec.~\ref{sec:exp:setup} followed by answering the three questions from Sec.~\ref{sec:exp:classic} to Sec.~\ref{sec:exp:generalize}. We also show additional experimental results and discussion reflecting condition number of the preconditioned system using various methods in Sec.~\ref{supp:cond}, comparing our proposed method with the multigrid approach in Sec.~\ref{supp:amg}, and learning physics simulation works in Sec.~\ref{supp:network-accu} in Appendix. Training setup can be found in Section~\ref{supp:setup}. in Appendix.

\subsection{Experiment Setup}\label{sec:exp:setup}

\paragraph{PDE Environments} This work studies solving the three representative linear second-order PDEs:
\begin{align}
& \text{Heat equation (parabolic)} \quad \pd{u}{t} - \alpha \pdtwo{u}{\xb}=0 \\
    & \text{Wave equation (hyperbolic)} \quad \pdtwo{u}{t} - c^2\pdtwo{u}{\xb}=0\\
    & \text{Poisson's equation (elliptic)} \quad \grad^2u = f,
\end{align}
where $\alpha$ is the diffusion coefficients, $c$ is the wave speed, and $f$ is the right hand side. Each of these PDEs is defined on a problem domain with a 2D triangle mesh and/or a 3D tetrahedron mesh for discretization purposes, as shown in Fig.~\ref{fig:env}). More details about the dataset generated for each environment can be found in Sec.~\ref{supp:sec:env} in Appendix. 

\begin{table*}[htb]
\resizebox{\textwidth}{!}{%
        \scriptsize
    \centering
    \begin{tabular}{lcccccccc}
    \toprule
Task & Method    & Precompute time 	$\downarrow$	& time (iter.)  $\downarrow$		& time (iter.) 	$\downarrow$	& time (iter.) 	$\downarrow$	& time (iter.) 	$\downarrow$	& time (iter.) 	$\downarrow$	& time (iter.)$\downarrow$ \\ 

& & (s) 		&  until 1e-2 		&  until 1e-4 		&  until 1e-6 		&  until 1e-8 		&  until 1e-10 		&  until 1e-12 \\
\midrule
\multirow{3}{*}{heat-2d}
&$\kappa$ losss (~\cite{sappl_arxiv})& 0.0236 &  0.232  (20) &  0.518  (68) &  0.652  (98) &  0.810  (118) &  1.023  (163) &  1.251  (180) \\

&Naive loss (Eqn (\ref{eqn:loss_A})) & 0.0218  & 0.177 (15) & 0.732 (181) & 1.187 (236) & 1.518 (296) & 1.885(296) & 2.284 (361)  \\
&Our loss (Eqn (\ref{eqn:loss_Ab})) & 0.0271 & \textbf{0.172} ({14}) & \textbf{0.420} ({57}) & \textbf{0.597} ({87}) & \textbf{0.733} ({110}) & \textbf{0.909} ({140}) & \textbf{1.056} ({165})\\
\midrule
\multirow{3}{*}{wave-2d}
&$\kappa$ losss (~\cite{sappl_arxiv}) &  0.0165 & 0.087 (0)  & 0.087 (0)  & 0.087 (0)  & 0.106 (3)  &  0.165 (12) & 0.220 (22) \\

&Naive loss (Eqn (\ref{eqn:loss_A})) & 0.0120  &  \textbf{0.076}  (0) &  \textbf{0.076} (0) &  \textbf{0.076} (0) &  0.146  (5) &  0.190  (21) &  0.351  (39) \\

&Our loss (Eqn (\ref{eqn:loss_Ab}))& 0.0147 & {0.081} (0) & {0.081} (0) & {0.081} (0) & \textbf{0.100} (3) & \textbf{0.156} (12) & \textbf{0.211} (21) \\

\midrule
\multirow{3}{*}{possion-2d}
&$\kappa$ losss (~\cite{sappl_arxiv}) &  0.0129 & 0.769(219) & 1.014 (291)  &  1.417(406) &  1.602 (471) &  1.990 (572) & 2.211 (662)  \\

&Naive loss (Eqn (\ref{eqn:loss_A})) & 0.0117 &  0.827 (231)  &  1.201 (319)  &  1.443 (413) &  1.634 (470)  & 1.942 (560)  & 2.256 (632)  \\

&Our loss (Eqn (\ref{eqn:loss_Ab})) & 0.0145 & \textbf{0.639} (175) &\textbf{0.818} (227) &\textbf{1.017} (286) &\textbf{1.118} (316) &\textbf{1.312} (374) &\textbf{1.510} (432)\\
\midrule
\multirow{2}{*}{poisson-3d}
&Naive loss (Eqn (\ref{eqn:loss_A})) & 0.407  & \textbf{2.936} (0)  & 3.241 (6)  &  \textbf{4.017} (21)  & 8.501 (46) & 12.373 (68)  &  15.719 (87) \\

&Ours (Eqn (\ref{eqn:loss_Ab})) & 0.413 & {3.010} (0) & \textbf{3.220} (2) & {4.815} (13) & \textbf{6.908} (28) & \textbf{8.749} (41) & \textbf{10.406} (53) \\
     \bottomrule
    \end{tabular}
    }
    \vspace{-2mm}
    \caption{Wall-clock time and iterations: our method with two different loss functions on {heat-2d}.}
    \label{tab:loss}
\end{table*}

\paragraph{Baseline Methods}
We compare with several general-purpose classical and learning-based preconditioners:
    
\begin{itemize}[leftmargin=*]
\setlength\itemsep{1mm}
    \item Jacobi preconditioner (Jacobi) uses the diagonal element of the original Matrix $\mathbf{A}$ as preconditioner $\mathbf{P}$, and its inverse can be easily computed by directly taking the inverse of the diagonal entries. 
    \item Gauss-Seidel preconditioner (Gauss-Seidel) is a factorization-based preconditioner. It constructs the upper $\mathbf{U}$ and lower $\mathbf{L}$ triangular matrices directly from the system matrix A, making $\mathbf{P}=\mathbf{L}+\mathbf{U}$. 
    \item Incomplete Cholesky preconditioner (IC) is a factorization-based preconditioner that is formed by the approximate triangular decomposition $\mathbf{P}=\mathbf{LL}^\top$. The numerical values $\mathbf{L}$ is sequentially computed from left-to-right to minimize $|| \mathbf{LL}^\top - \mathbf{A} ||_2 $. 
    \item Incomplete Cholesky with two levels of fill-in (IC(2)) is a variant of the standard IC preconditioner that improves accuracy. It uses the second level of fill-in to capture additional non-zero entries in the factorization. It first uses IC to factorize $\mathbf{A}$ into $\mathbf{LL}^{\top}$. Then, nonzero entries of L are used to construct a new sparse matrix $\mathbf{A}_2$, which is factorized again using IC, resulting in a new sparse matrix $\mathbf{L}_2$. The final preconditioner is  $\mathbf{L}_2\mathbf{L}_2^\top$. IC(2) is more accurate than IC as it captures more of the structure of $\mathbf{A}$. However, it is also more computationally expensive to construct, requiring two rounds of factorizations. 
    \item Learning-based preconditioners trained by directly minimizing the condition number \citep{sappl_arxiv}.
\end{itemize}

\begin{table*}[htb]
\resizebox{\textwidth}{!}{%
        \scriptsize
\begin{tabular}{lcccccccc}

\toprule
Task & Method    & Precompute time 	$\downarrow$	& time (iter.)  	$\downarrow$	& time (iter.) 	$\downarrow$	& time (iter.) 	$\downarrow$	& time (iter.) 	$\downarrow$	& time (iter.) 	$\downarrow$	& time (iter.) $\downarrow$ \\ &

& (s) 		&  until 1e-2 		&  until 1e-4 		&  until 1e-6 		&  until 1e-8 		&  until 1e-10 		&  until 1e-12 \\

\midrule
\multirow{5}{*}{test}

 &Jacobi & \textbf{0.0001} & 0.811 (230) & 0.983 (281) & 1.341 (388) & 1.629 (474) & 1.832 (535) & 2.176 (640) \\
&Gauss-Seidel & 0.0065 &  0.576  (157) &  0.757  (211) &  1.076  (305) &  1.333  (382) &  1.6  (461) &  1.897  (550) \\
 &IC & 0.6952 &  1.079  ({100}) &  1.157  ({123}) &  1.275  ({156}) &  1.437  ({203}) &  1.525  ({229}) &  1.648  ({264})  \\

&IC(2) & 1.1123 &  1.354  (60) &  1.396  (72) &  1.491  (98) &  1.556  (115) &  1.607  (129) &  1.707  (156) \\

&Ours & 0.0138 & \textbf{0.513} (137) & \textbf{0.615} (167) & \textbf{0.814} (226) & \textbf{0.987} (277) & \textbf{1.109 }(313) & \textbf{1.297} (369)  \\

\midrule
\multirow{5}{*}{test-$\sigma$} 

&Jacobi & \textbf{0.0001 }&  0.904  (253) &  1.1  (310) &  1.505  (430) &  1.756  (505) &  2.041  (589) &  2.414  (702)  \\ 

&Gauss-Seidel & 0.007 &  0.638  (175) &  0.814  (227) &  1.195  (340) &  1.425  (408) &  1.771  (512) &  2.046  (594) \\

&IC & 0.7093 &  1.104  ({110}) &  1.178  ({132}) &  1.338  ({180}) &  1.461  ({217}) &  1.551  ({244}) &  1.715  ({293}) \\

&IC(2) & 1.1305 &  1.362  (57) &  1.403  (69) &  1.497  (94) &  1.562  (112) &  1.611  (125) &  1.708  (151) \\

&Ours & 0.0139 & \textbf{0.603} (165) & \textbf{0.72} (200) & \textbf{0.986} (280) & \textbf{1.14} (326) & \textbf{1.292} (372) & \textbf{1.557} (452)  \\

\midrule
\multirow{5}{*}{test-$3\sigma$} 

&Jacobi &\textbf{ 0.0001} &  0.949  (262) &  1.139  (317) &  1.566  (441) &  1.828  (519) &  2.123  (606) &  2.485  (714) \\

&Gauss-Seidel & 0.0064 &  0.588  (161) &  0.749  (209) &  1.133  (323) &  1.334  (383) &  1.694  (491) &  1.942  (565) \\

&IC& 0.6992 &  1.114  ({111}) &  1.198  ({137}) &  1.365  ({187}) &  1.484  ({222}) &  1.578  ({251}) &  1.743  ({301}) \\

&IC(2) & 1.1126 &  1.354  (60) &  1.396  (72) &  1.491  (98) &  1.556  (115) &  1.608  (129) &  1.707  (156) \\

&Ours & 0.0141 & \textbf{0.56} (153) & \textbf{0.672} (186) & \textbf{0.905} (256) & \textbf{1.051} (300) & \textbf{1.201} (345) & \textbf{1.433} (415) \\

\midrule
\multirow{5}{*}{test-$5\sigma$} &

Jacobi &\textbf{ 0.0001} &  1.082  (310) &  1.455  (421) &  1.761  (512) &  2.06  (603) &  2.531  (745) &  2.763  (816)  \\
&Gauss-Seidel & 0.0065 &  \textbf{0.66}  (182) &  \textbf{0.836}  (234) &  1.218  (347) &  1.455  (417) &  1.805  (522) &  2.08  (604) \\

&IC & 0.6939 &  1.188  ({135}) &  1.311  ({171}) &  1.472  ({219}) &  1.564  ({246}) &  1.783  ({311}) &  1.912  ({349}) \\

&IC(2) & 1.1019 &  1.317  (53) &  1.359  (65) &  1.444  (88) &  1.513  (107) &  \textbf{1.561}  (119) &  \textbf{1.657}  (146) \\

&Ours & 0.0136 & {0.728} (203) & {0.97} (275) & \textbf{1.165} (333) & \textbf{1.334} (384) & {1.683} (488) & {1.817} (528) \\

\bottomrule
\end{tabular}
}
\vspace{-2mm}
    \caption{Generalization to different physics parameters. We test on test sets with increasing deviation from the training distribution. $\sigma$ stands for the standard deviation of training set, test-$\sigma$, test-$3\sigma$, and test-$5\sigma$ means test sets that are of 1, 3, and 5 std. dev away from the training distribution, respectively. }
    \label{tab:generalize}
    \vspace{-2mm}
\end{table*}

We also include a detailed discussion between our method and the multigrid preconditioner method in Appendix~\ref{supp:amg}. More details about baselines can be found in Appendix.

\paragraph{Evaluation Metrics}

We quantify the performance by comparing the total wall-clock time spent for each preconditioner to reach desired accuracy levels.

To ensure a fair comparison between all methods, we summarize the performance of PCG solvers not in a single number but in the following values:
(a) the time spent on precomputing the preconditioner for the given $\mathbf{A}$ and $\mathbf{b}$;
(b) the number of iterations for CG solver to converge, and (c) the total time (including the precomputing time) to reach different precision thresholds.

As an additional reference metric, we also show the condition number that reflects the speed up during the CG solving stage in Appendix Sec.~\ref{supp:cond}. 

\begin{table*}[htb]
\centering
\resizebox{\textwidth}{!}{%
        \scriptsize
\begin{tabular}{cccccccc}
\toprule
Method    & Precompute time $\downarrow$		& Time (iter.)  	$\downarrow$	& Time (iter.) 	$\downarrow$	& Time (iter.) 	$\downarrow$	& Time (iter.) 	$\downarrow$	& Time (iter.) 	$\downarrow$	& Time (iter.) $\downarrow$ \\

&  		&  until 1e-2 		&  until 1e-4 		&  until 1e-6 		&  until 1e-8 		&  until 1e-10 		&  until 1e-12 \\

\midrule

Jacobi & \textbf{0.0002} &  {2.314}  (6) &  {4.634}  (22) &  8.395  (48) &  12.002  (72) &  15.381  (95) &  18.332  (115) \\

Gauss-Seidel & {0.3167} &  \textbf{0.886}  (0) &  \textbf{1.728}  (8) &  \textbf{5.54}  (29) &  \textbf{8.768}  (48) &  11.83  (65) &  14.847  (82) \\
IC & 8.9818 &  9.686  (2) &  10.549  (12) &  11.352  (21) &  12.097  (29) &  12.721  (36) &  13.524  (45) \\
IC(2) & 14.0376 &  14.688  (2) &  15.475  (10) &  16.072  (16) &  16.706  (23) &  17.269  (29) &  18.008  (37) \\
Ours & 0.4206 & {3.591} (4) & {4.97} (14) & {7.005} (29) & {8.978} (43) & \textbf{10.721} (55) & \textbf{12.412} (68)  \\

\bottomrule
\end{tabular}
}
\vspace{-2mm}
\caption{Our approach generalizes to unseen armadillo mesh in the {poisson-3d} environment.}
\label{tab:mesh}
\vspace{-2.5mm}

\end{table*}

\subsection{Comparison with Classic Preconditioners}\label{sec:exp:classic}
We compare our approach with the general-purpose preconditioners described above. Table~\ref{tab:pcg} summarizes the time cost and iteration numbers of PCG solvers using different preconditioners up to various convergence thresholds.

Jacobi conducts the simple diagonal preconditioning, so it has little precomputation overhead. However, the quality of the preconditioner is mediocre, and PCG takes many iterations to converge. 
Similarly, the Gauss-Seidel preconditioner directly forms the triangular matrices by taking the entry values of system matrix $\mathbf{A}$, and thus the computational overhead is low as compared to IC or IC(2) preconditioners. However, it is limited in speeding up solvers because of its coarse approximation to the system matrix $\mathbf{A}$.

IC and IC(2) speeds up CG solver significantly. However, its precomputation process is sequential and therefore expensive to compute. This comparison reflects the derivation complexity and approximation accuracy trade-off between existing numerical preconditioners.  By contrast, our approach features an easily parallelizable precomputation stage (like Jacobi and Gauss-Seidel) and produces a preconditioner with a quality close to IC. Therefore, in terms of total computing time, our approach outperforms the existing general-purpose numerical approaches across a wide range of precision thresholds.

We can also see from this experiment that our proposed method is especially beneficial for large-scale problems in Poisson-3d with a matrix size of $23300 \times 23300$. We outperform baseline methods by a large margin. This reflects the clear advantage of our approach's parallelizability as opposed to sequential approaches such as IC or IC(2). 

\subsection{Comparison of Loss function}\label{sec:exp:loss}
To highlight the value of our loss function targeting data distributions on the training set, we train the network with the loss function that reduces the condition number $\kappa$ as proposed in~\cite{sappl_arxiv}, the loss function that focuses only on the system matrix $\mathbf{A}$ as shown in Eqn (\ref{eqn:loss_A}), and our proposed loss function that introduces inductive bias as shown in Eqn (\ref{eqn:loss_Ab}). We show the comparisons across all four of our experimental settings, and the results are shown in Table~\ref{tab:loss}. Compared to Eqn (\ref{eqn:loss_A}), We observe that the preconditioner trained with Eqn (\ref{eqn:loss_Ab}) converges in fewer iterations than Eqn (\ref{eqn:loss_A}). 

Compared with the $\kappa$ loss function, our proposed data-driven loss function (\ref{eqn:loss_Ab}) shows a more stable convergence across various different second-order linear PDEs.  We also observe that $\kappa$ loss~\cite{sappl_arxiv} works on par with our method on the wave-2d setting, but it does not work in other settings, e.g., heat-2d, poisson-2d. Additionally, the method using condition number as the loss energy is not computationally efficient. Computing condition number scales cubically with problem sizes. Our poisson-3d setting uses a mesh of size 23,300 nodes. We train the $\kappa$ loss method for five days (120 hours), but it does not converge, and it is only trained through less than five percent of the training set using the same hardware setup. CG algorithm does not converge when testing on these non-convergent results. Empirically, we found that using condition number as loss energy is computationally infeasible with problems of more than 10,000 nodes. As such, we conclude that enforcing data distribution dependence during training allows us to achieve better in-distribution inference during test time.

\subsection{Generalization}\label{sec:exp:generalize}

\textbf{Physics parameters.} First, we consider generalizing the PDEs on their physics parameters, which govern system $\mathbf{A}$. We use {poisson-2d} as an example. Our method is trained on a fixed density distribution between $[0.001, 0.005]$, and we test the performance of our method on test distributions that gradually deviates from the training distribution. Results are reflected in Table~\ref{tab:generalize}, with growing deviation from top rows to bottom rows. 

Since changing physics parameters does not affect the matrix sparsity, the pre-computation time remains largely unchanged across different data distributions (see Table~\ref{tab:generalize} Column 1). We can see that even on the challenging out-of-distribution datasets, our approach still maintains reasonable performance. We achieve high precision while using the least total time to maintain better or comparable performance to existing numerical approaches. We also observe that the performance of our method degrades as the domain gap between the training and test distribution grows.

\textbf{Geometry.} We test the generalizability of our model on different problem domains $\Omega$, represented by different mesh models in our setting. We train our network on the connector shape on one mesh and test its performance on unseen meshes. As shown in the {poisson-3d} environment, we train the network preconditioner on a ``connector mesh'' consisting of $23,300$ nodes. We then test the trained network model on a new mesh ``armadillo'', which consists of $18,181$ nodes. The geometry of the two domains differs significantly as shown in Fig.~\ref{fig:env}.
We also report the time and the iteration cost of our approach and classic preconditioners on {poisson-3d} when testing on the unseen armadillo mesh. The results are shown in Table~\ref{tab:mesh}. We notice that our method can generalize to the unseen mesh. Even in the challenging case of solving PDE on an unseen test mesh, we are still able to converge to high precision levels faster compared to existing numerical approaches. We show additional experiments comparing our method with pure learning-based methods~\cite{pfaff2020learning} on generalizability in Appendix~\ref{supp:exp:network}. Since our approach is embedded inside the PCG framework, we significantly outperform these pure learning approaches.

\section{Conclusions and Future Work}\label{sec:conclusion}
This work presents a generic learning-based framework for estimating preconditioners in the context of conjugate gradient PDE solvers. 
Our key observation is that the preconditioner for classic iterative solvers does not require exact precision and is an ideal candidate for neural network approximation. Our proposed method approximates the preconditioner with a graph neural network and embeds this preconditioner into a classic iterative conjugate gradient solver. 
Compared to classic preconditioners, our approach is faster while achieving the same accuracy.

Currently, our approach is limited to linear systems of equations $\mathbf{Ax}=\mathbf{b}$. Our parallelizability is bounded by hardware setup, i.e., GPU memory. Future work may consider extending to dynamic sparsity patterns for larger systems and to more complex PDEs, such as the elastodynamics equations shown in prior end-to-end ML approaches \citep{sanchez2020learning}. 

\section{Acknowledgement}
The work is supported by the MIT Robert Shillman Fellowship.

\bibliography{iclr23}
\bibliographystyle{icml2023}

\clearpage
\appendix
\section{Appendix}
The Appendix section includes the following:
\begin{itemize}[leftmargin=*]
    \item Algorithm Description for PCG
    \item Details of the Experiment Environment
    \item \revised{Technical Details of our Method}
    \item
    \item \revised{Training Setup and Convergence Time}
    \item
    \item \revised{Additional Experiments on Large Matrices}
    \item
    \item \revised{Comparison with Start Guess $x_0$ Prediction}
    \item 
    \item Condition Number Comparison
    \item Comparison and Discussion on Multigrid Preconditioners
    \item Generalizability Comparison with Learning Physics Simulation Works
    \item Error Accumulation Comparison with Learning Physics Simulation Works
\end{itemize}

\subsection{PCG Algorithm}
\label{supp:sec:pcg}
We describe the Preconditioned Conjugate Gradient Algorithm Here:

\begin{algorithm}[htb]
\caption{PCG}
\label{alg:pcg}
\begin{algorithmic}
\REQUIRE{System Matrix $A$, right-hand-side vector $b$, initial guess $x_0$, Precondtioner $P$}
\STATE ${\rbd}_0=\bb-A\xb_0$
\STATE Solve $P\zb_0=\rbd_0$
\STATE $\pb_1 = \zb_0$
\STATE $w = A\pb_1$
\STATE $\alpha_1 = \rbd_0^T\zb_0 / (\pb_1^T\wb)$
\STATE $\xb_1 = \xb_0 + \alpha_1\pb_1$
\STATE $\rbd_1 = \rbd_0 - \alpha_1\wb$
\STATE $k=1$
\WHILE{$\|\rbd_k\|_2 > \epsilon$}
\STATE Solve $P\zb_k=\rbd_k$
\STATE $\beta_k = \rbd_k^T\zb_k / (\rbd_{k-1}^T\zb_{k-1})$
\STATE $\pb_{k+1} = \zb_k + \beta_k\pb_k$
\STATE $w = A\pb_{k+1}$
\STATE $\alpha_{k+1} = \rbd_k^T\zb_k / (\pb_{k+1}^T\wb)$
\STATE $\xb_{k+1} = \xb_k + \alpha_{k+1}\pb_{k+1}$
\STATE $\rbd_{k+1} = \rbd_k - \alpha_{k+1}\wb$
\STATE $k=k+1$
\ENDWHILE
\STATE \textbf{Return} $\xb = \xb_k$
\end{algorithmic}
\end{algorithm}

\subsection{Environment Details}
\label{supp:sec:env}
Our dataset is constructed by simulating trajectories of each PDE on a given mesh domain with various initial conditions and boundary conditions. Each trajectory has about 20-100 time steps depending on the equation and initial condition. These trajectories are split into training and test set. ${\mathbf{A}_i, \mathbf{x}_i, \mathbf{b}_i}$ is one single time step in the trajectory. The training set is of size 3000, meaning that it contains 3000 $\mathbf{A}_i \cdot \mathbf{x}_i = \mathbf{b}_i$ tuples, and the test set contains 200 instances.

\begin{table}
    \centering
    \begin{scriptsize}
    \begin{tabular}{c|c|c|c}
    \toprule
        PDE &  Number of Nodes & Number of Elements & Boundary Condition\\ 
    \midrule
        heat-2d & 7454 & 14351  & varying \\
        wave-2d & 4852 & 8839 & varying \\
        poisson-2d&  3167 & 6117 & varying \\
        poisson-3d &  23300 & 129981 &  fixed \\
        poisson-3d &  18181 & 97476 &  fixed \\
    \bottomrule
    \end{tabular}
    \end{scriptsize}
    \caption{Environment setup for experiments.}
    \label{tab:environemnt}
\end{table}

\begin{figure*}[htb]
    \centering
    \includegraphics[width=\textwidth]{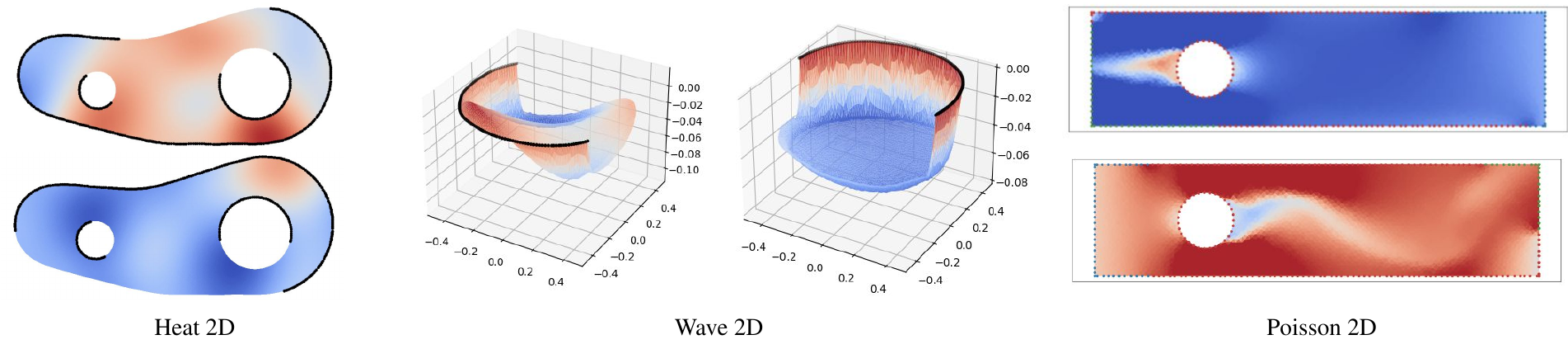}
    \caption{\small{Varying boundary conditions: For {heat-2d} and {wave-2d}, the black vertices represent the Dirichlet boundary. The {poisson-2d} example shows the red vertices are the obstacle boundary, the blue vertices are the influx boundary, and green vertices are the Dirichlet boundary.}}
    \label{supp:fig:supp-data}
\end{figure*}

Table~\ref{tab:environemnt} lists the environment details for the experiment section.
Figure~\ref{supp:fig:supp-data} shows examples demonstrating varying boundary conditions. For {heat-2d} and {wave-2d}, varying lengths and positions of mesh geometric boundary nodes are selected as Dirichlet boundaries. 
For {poisson-2d} equation, we use the inviscid-Euler fluid equation as a demonstration. All solvers are only responsible for solving the pressure that makes the velocity field incompressible, which is a Poisson equation. The advection and external force steps are then applied to generate the data visualization. For {poisson-2d}, two sets of varying length and position of mesh geometric outer border boundary nodes are selected as influx and Dirichlet boundary. The remaining mesh geometric border nodes, including the remaining outer border and all nodes in the inner border, are obstacle boundaries. 

For experiment~\ref{sec:exp:generalize} across different physics parameters, we consider the same mesh domain used in other {poisson-2d} experiments with a mesh size of 3167 and element size of 6117. We train on training sets with density distribution from [0.001, 0.005], and our test environment \textit{test-$1\sigma$} is of density 0.006. Test environment \textit{test-$3\sigma$} is of density 0.008. Test environment \textit{test-$5\sigma$} is of density 0.01.

\subsection{Technical Details and Justification} 
\label{supp:sec:ours}
\revised{
Encoders operate on graph nodes and edges. Graph Node Encoder is a $l$ layer MLP with $h$ hidden dimensions that takes each graph node input (rhs vector $\mathbf{b}$) to a 16-dimensional feature vector. Graph Edge Encoder is an MLP also with $l$ layers and $h$ hidden dimensions. It operates on graph edges input (matrix A) to a 16-dimensional latent feature. }

\revised{
Message passing of $n_{mp}$ iterations is conducted where the neighboring nodes are updated through the connected edges. $v_{i, t+1} = f_{mp, v} ( v_{i, t}, \sum_j e_{i, j, t} v_{j, t})$ where $f_{mp, v}$ is implemented as an MLP with $l_{mp}$ layers and $h_{mp}$ hidden dimensions. The Edge features are also updated by combining the updated information of the two connected nodes through the message passing layers such that $e_{i, j, t+1} = f_{mp, e} ( e_{i, j, t} ,  v_{i, t+1}, v_{j, t+1} )$,  where $f_{mp, e}$ is also implemented as an MLP with $l_{mp}$ layers and a $h_{mp}$ hidden dimensions. We then update the two-way edges $e_{i, j, t+1}$ and $e_{j, i, t+1}$.}

 \revised{ The decoder converts the updated edge features $e_{i, j}$ into a single real number $L_{i, j}. L_{i, j} = f_{dec} (e_{i, j})$, where $f_{dec}$ a $l$ layer MLP with $h$ hidden dimensions. To ensure the decoded $\mathbf{L}$ is a lower triangular matrix, we average the value on the pair of bi-directional edges $e_{i,j} = \frac{1}{2}(e_{i, j} + e_{j, i})$ and store the value on the corresponding lower-triangular indices $\mathbf{L}_{i, j | i \geq j }$. The output of the GNN is $\mathbf{L}$. ReLU activation is used in all MLPs. }

\begin{table}[htb]
    \centering
    \begin{scriptsize}
    \begin{tabular}{c|cccc}
    \toprule
        Env &  Heat-2D & Wave-2D & Poisson-2D & Poisson-3D \\
        \midrule
        $l$ & 1 & 1 & 2 & 2 \\
         $h$ & 16 & 16 & 16 & 16  \\
         $n_{mp}$ &  5 & 5 & 5 & 3 \\
         $l_{mp}$ &  1 & 1 & 2 & 2   \\
         $h_{mp}$ &  16 & 16 & 16 & 16 \\
        \bottomrule
    \end{tabular}
    \caption{{GNN architecture hyper-parameter.}}
    \label{supp:gnn-arch}
    \end{scriptsize}
\end{table}

We follow the diagonal decomposition $\mathbf{L}\mathbf{D}\mathbf{L}^{\top}$ as a way of lower triangular decomposition for the original system $\mathbf{K}$. It is easy to see that this diagonal decomposition is equivalent to lower triangular decomposition. 
\begin{align}
\mathbf{K} &= \mathbf{L}_{\theta} \mathbf{L}_{\theta}^{\top} \\
& =  \mathbf{L}_{\theta'} \sqrt{\mathbf{D}} \sqrt{\mathbf{D}} \mathbf{L}_{\theta'} \\
& = \mathbf{L}_{\theta'} \mathbf{D} \mathbf{L}_{\theta'}
\end{align}
The diagonal decomposition $\mathbf{L}\mathbf{D}\mathbf{L}^{\top}$ has several advantages, similar to lower triangular decomposition $\mathbf{L}\mathbf{L}^{\top}$, it is easy to invert, and guarantees symmetry. Additionally, we enforce the diagonal element $\mathbf{D}$ to be the diagonal elements of the original system $\mathbf{K}$. This way, we enforce the value and gradient range for the lower triangular matrix $\mathbf{L}_{\theta'}$ to ensure the positive definiteness of the learned decomposition.

\subsection{Training Time and Training Setup Details}
\label{supp:setup}


\revised{\textbf{Training Setup. }  
All experiments are conducted using the same hardware setup equipped with 64-core AMD CPUs and an NVIDIA RTX-A8000 GPU. We use Adam optimizer with the initial learning rate set to 1e-3. We use a batch size of 16 for all Heat-2d, Poisson-2d, and Wave-2d experimental environments. We use a batch size of 8 for the Poisson-3d environment. All learning-based methods are written using the Pytorch~\cite{NEURIPS2019_9015} Framework with the Pytorch-Geometric~\cite{pytorch-geometric} Package and CUDA11.6. We use the same set of GNN architecture hyperparameters as described in Sec.~\ref{supp:gnn-arch} for all four experimental environments. 

All learning-based methods are trained to full convergence unless otherwise specified. All learning-based and traditional baselines use parallelized tensorized GPU implementations unless otherwise specified. All experimental data are reported by averaging the test sets of size 200. Each experiment is run 12 times, and the reported time averages the fastest five runs. 
We trained our proposed method for 5 hours on all experimental environments and observed full convergence. All learning baseline methods discussed in Sec. A~\ref{sec:exp:loss} are trained for 72 hours and observed convergence except for the Poisson-3d environment, which is trained for 120 hours and does not converge.  }

\begin{figure}[h]
     \centering
     \includegraphics[width=0.8\linewidth]{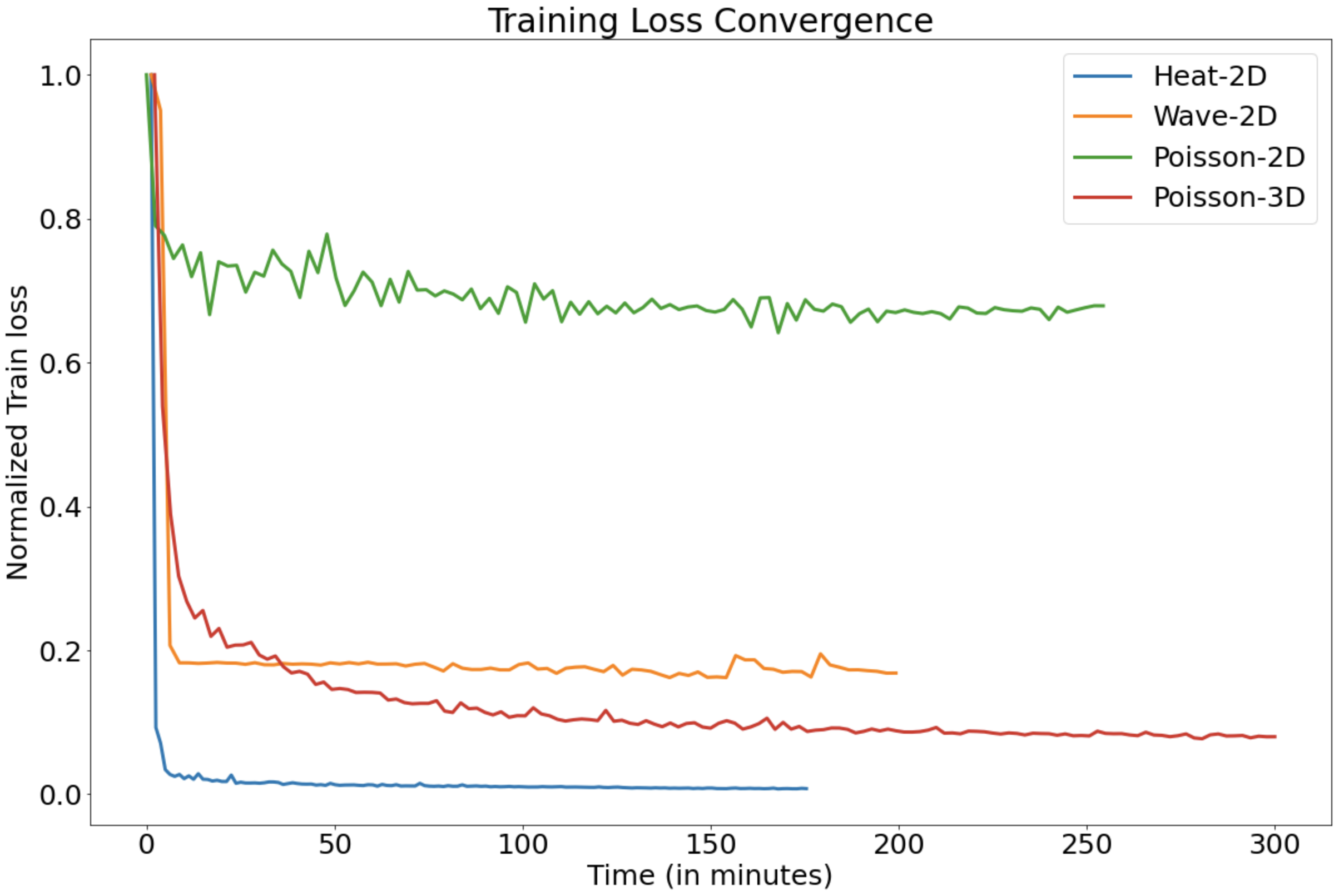}
     \caption{Training Convergence (in normalized loss value) }
    \label{supp:fig:convergence-Comparison}
\end{figure}

\revised{\textbf{Training Time. }  Figure~\ref{supp:fig:convergence-Comparison} reflects the convergence speed for training each PDE environment. We train each PDE environment for 5-6 hours. We observe that training for all environments converges within one hour of training. Heat-2d environments converge within 5 minutes, and Wave-2d environments converge within 8 minutes. The Poisson-2D and Poisson-3d environments converge within about 30 and 60 minutes of training, respectively. It is noted that we want the wall-clock time value reported in all tables to be one single simulation time step. A full simulation trajectory is normally composed of thousands of or even millions of such time steps. Therefore, the small time difference for each simulation time step can quickly add up to pay off the training cost.  }

\subsection{Additional Experiments on Large Matrices}
\label{supp:large}
\revised{We provide several additional large-scale examples ($>10,000$ mesh nodes/matrix dimension) to show the advantages of our method. The experimental results are shown in Table~\ref{tab:large-exp}. We can see from the table that our proposed method surpasses the classical baseline methods by a large margin, especially on the low-accuracy requirement tasks ($1e-2 \sim 1e-6$). }

\begin{table*}[htb]
\resizebox{\textwidth}{!}{%
        \scriptsize
\centering
\begin{tabular}{lcccccccc}

\toprule
Task &  Mesh Shape & Method  & Precompute time 	$\downarrow$	& time (iter.) $\downarrow$ 		& time (iter.) 	$\downarrow$	& time (iter.) 	$\downarrow$	& time (iter.) $\downarrow$		& time (iter.) 	$\downarrow$	 \\ 
 problem size  &  & & (s) 		&  until 1e-2 		&  until 1e-4 		&  until 1e-6 		&  until 1e-8 		&  until 1e-10 		\\

\midrule
\multirow{5}{1.2 cm}{heat-2d (13260) } & \multirow{5}{*}{  \begin{minipage}{1.3cm}
      \includegraphics[width=1.5cm, height=4mm]{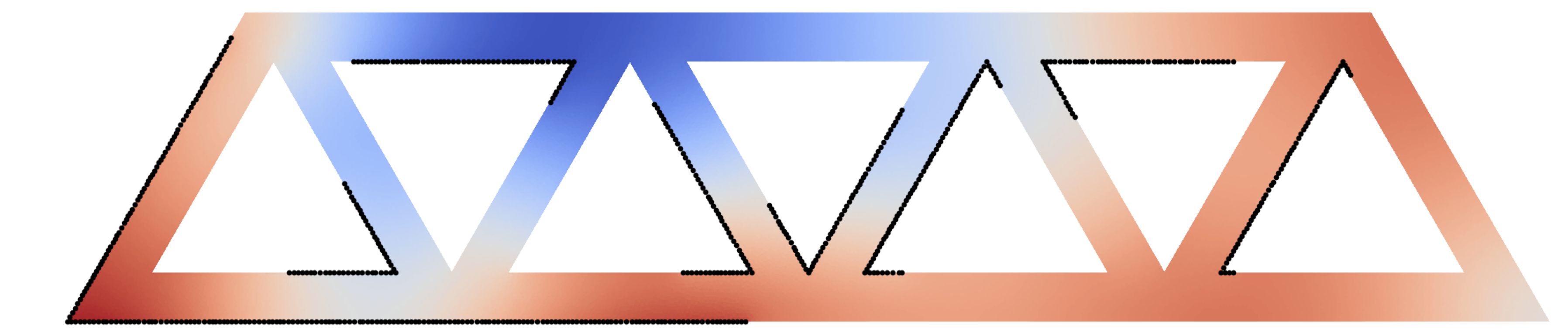}
    \end{minipage}} &

Jacobi & \textbf{0.0001} &  2.036  (33) &  6.93  (136) &  10.403  (210) &  13.427  (274) &  16.931  (348) \\
&&Gauss-Seidel & 0.0688 &  1.931  (23) &  6.385  (95) &  9.545  (147) &  12.297  (191) &  15.486  (244) \\
&&IC & 3.4732 & 4.911 ({13}) & 6.947 ({56}) & 8.428 ({87}) & 9.541 ({110}) & 10.852 ({138})  \\
&&IC (2) & 5.8823 &  7.332  (9) &  8.961  (39) &  10.146  (61) &  11.036  (77) &  12.085  (96)  \\
&&Ours & 0.0478 & \textbf{1.607} (17) & \textbf{4.088} (73) & \textbf{5.222} (115) & \textbf{7.084} (146) & \textbf{8.057} (184) \\

\midrule
\multirow{5}{1.2 cm}{possion-2d (13436)} & \multirow{5}{*}{ \begin{minipage}{1.3cm}
      \includegraphics[width=1.3cm, height=1.0cm]{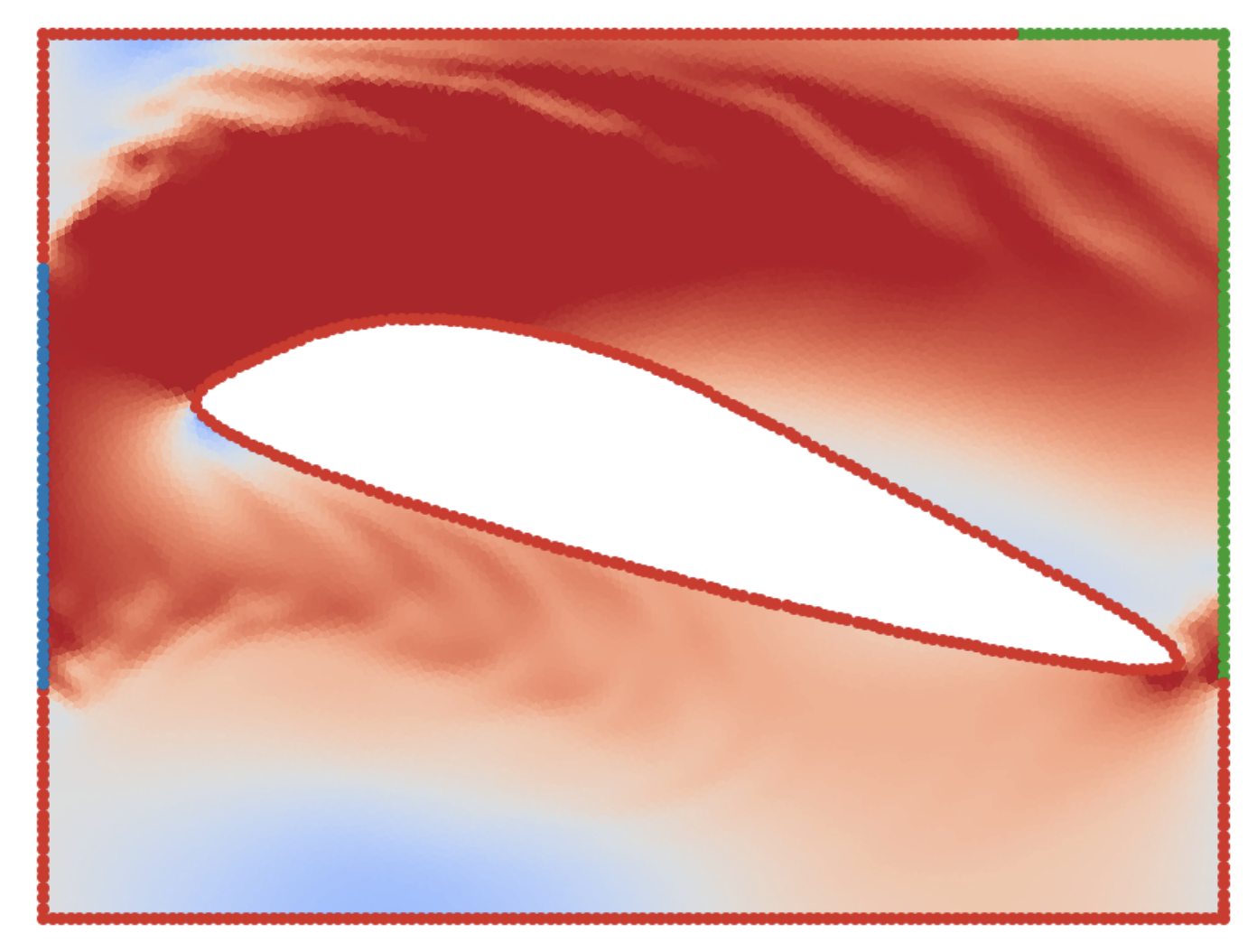}
    \end{minipage}}

&Jacobi & \textbf{0.0001} &  18.929  (496) &  27.988  (733) &  35.38  (927) &  45.646  (1194) &  52.419  (1371) \\
&&Gauss-Seidel & 0.0594 &  12.577  (327) &  18.558  (484) &  23.454  (612) &  30.192  (788) &  34.659  (905)  \\
&&IC & 3.7626 & 11.859 ({212}) & 13.776 ({262}) & 17.958 ({371}) & 20.042 ({426}) & 23.042 ({504})  \\
&&IC (2) & 6.6974 &  13.008  (144) &  14.514  (178) &  17.876  (252) &  19.241  (290) &  21.941  (343)  \\
&&Ours & 0.0427 & \textbf{9.607} (256) & \textbf{10.115} (306) & \textbf{15.07} (443) & \textbf{11.477} (524) & \textbf{16.924} (576)  \\

\bottomrule

\end{tabular}
}
\vspace{-3mm}
\caption{Comparison between preconditioners for PCG on large examples. We report the precompute time, total time ( ICl. precompute time) for each precision level, and the corresponding PCG iterations (in parenthesis). The best value in each category is in bold. $\downarrow$: the lower the better.}
\label{tab:large-exp}
\end{table*}

\subsection{Comparison with Starting Guess $\hat{x}_0$ prediction}
\label{supp:x-compare}
\revised{Here, we provide an additional experiment using both our predicted preconditioner $\mathbf{P}$ and the initial value of $\hat{\mathbf{x}}_0$ for the PCG algorithm. $\hat{\mathbf{x}}_0$ is obtained by regressing the decoded graph node values to the $\mathbf{x}$ in $\mathbf{Ax}=\mathbf{b}$. The starting guess $\hat{\mathbf{x}}_0$ can be obtained with no additional computation time cost. We compare with existing classical preconditioning methods as mentioned in Section~\ref{sec:exp:classic}. Experimental results are shown in Table~\ref{tab:comp-x}. }

\begin{table*}[htb]
\resizebox{\textwidth}{!}{%
        \scriptsize
\centering
\begin{tabular}{lcccccccc}

\toprule
Task & Method    & Precompute time 	$\downarrow$	& time (iter.) $\downarrow$ 		& time (iter.) 	$\downarrow$	& time (iter.) 	$\downarrow$	& time (iter.) $\downarrow$		& time (iter.) 	$\downarrow$	& time (iter.)$\downarrow$ \\ &

& (s) 		&  until 1e-2 		&  until 1e-4 		&  until 1e-6 		&  until 1e-8 		&  until 1e-10 		&  until 1e-12 \\

\midrule
\multirow{6}{*}{heat-2d} &

 Jacobi & \textbf{0.0001} &  0.657  (32) &  2.188  (132) &  3.263  (202) &  4.105  (257) &  5.269  (333) &  6.255  (398) \\
&Gauss-Seidel & 0.0071 &  0.645  (27) &  1.656  (98) &  2.339  (146) &  2.995  (193) &  3.771  (247) &  4.402  (292) \\

& IC& 1.5453 &  1.954  ({12}) &  2.612  ({54}) &  3.061  ({83}) &  3.409  ({105}) &  3.832  ({133}) &  4.271  ({161}) \\

&IC(2) & 2.3591 &  2.624  (11) &  3.094  (40) &  3.399  (60) &  3.651  (76) &  3.965  (96) &  4.260  (115) \\
&Ours & 0.0251 & {0.490} (17) & {1.284} (71) & {1.856} (110) & {2.300} (140) & {2.831} (177) & {3.377} (214) \\

& Ours with $\hat{x}_0$ & 0.0237 & \textbf{0.413} (14) & \textbf{1.24} (72) & \textbf{1.782} (110) & \textbf{2.219} (141) & \textbf{2.766} (180) & \textbf{3.284} (216) \\

\midrule
\multirow{6}{*}{wave-2d}
& Jacobi & \textbf{0.0001} &  0.141  (0) &  0.141  (0) &  0.141  (0) &  0.176  (6) &  0.295  (26) &  0.417  (46) \\

&Gauss-Seidel & 0.0079 &  0.089  (0) &  0.089  (0) &  0.09  (0) &  0.1  (2) &  0.16  (11) &  0.232  (22) \\

& IC& 0.7679 &  0.885  ({0}) &  0.885  ({0}) &  0.885  ({0}) &  0.904  ({3}) &  0.953  ({11}) &  1.007  ({20}) \\

&IC(2) & 1.1831 &  1.226  (0) &  1.226  (0) &  1.226  (0) &  1.266  (5) &  1.326  (12) &  1.385  (18) \\
&Ours & 0.0147 & {0.081} (0) & {0.081} (0) & {0.081} (0) & {0.100} (3) & {0.156} (12) & {0.211} (21) \\

&Ours with $\hat{x}_0$ & 0.0143 & \textbf{0.079} (0) & \textbf{0.079} (0) & \textbf{0.079} (0) & \textbf{0.097} (3) & \textbf{0.147} (11) & \textbf{0.201} (20) \\

\midrule
\multirow{6}{*}{possion-2d} &

 Jacobi & \textbf{0.0001} & 0.980 (275) &1.231 (348) &1.572 (448) &1.822 (522) &2.119 (611) &2.405 (697)\\

&Gauss-Seidel & 0.0071 &  0.699  (194) &  0.964  (273) &  1.26  (361) &  1.518  (438) &  1.807  (525) &  2.099  (613) \\

& IC& 0.7093 &  1.188  ({135}) &  1.309  ({171}) &  1.468  ({219}) &  1.559  ({246}) &  1.774  ({311}) &  1.900  ({349}) \\

&IC(2) & 1.205 &  1.308  (60) &  1.439  (74) &  1.543  (100) &  1.604  (115) &  1.664  (131) &  1.747  (151) \\
&Ours & 0.0145 & {0.639} (175) &{0.818} (227) &{1.017} (286) &{1.118} (316) &{1.312} (374) &{1.510} (432)\\

&Ours with $\hat{x}_0$ & 0.0137 & \textbf{0.588} (167) & \textbf{0.701} (203) & \textbf{0.952} (282) & \textbf{1.083} (322) & \textbf{1.276} (383) & \textbf{1.494} (451) \\
\midrule
\multirow{6}{*}{possion-3d}

& Jacobi & \textbf{0.0002} & \textbf{1.526} (0) & \textbf{2.693} (7) & 5.496 (25) & 9.552 (50) & 13.636 (76) & 17.080 (97) \\

&Gauss-Seidel & 0.3381 &  5.824  (0) &  6.775  (6) &  9.074  (19) &  12.305  (38) &  15.454  (56) &  18.333  (72) \\
&IC & 9.6878 &  10.668  (1) &  11.353  (6) &  12.592  (15) &  13.826  (23) &  14.954  (31) &  15.812  (37) \\
&IC(2) & 17.138 &  18.083  (1) &  18.661  (5) &  19.667  (11) &  20.599  (17) &  21.704  (24) &  22.560  (30) \\
&Ours & 0.4137 & {3.010} (0) & {3.220} (2) & {4.815} (13) & {6.908} (28) & {8.749} (41) & {10.406} (53) \\
&Ours with $\hat{x}_0$ & 0.4068 & {2.997} (0) & {3.209} (2) & \textbf{4.803} (12) & \textbf{6.882} (27) & \textbf{7.915} (40) & \textbf{10.020} (51) \\
\bottomrule

\end{tabular}
}
\vspace{-3mm}
\caption{Comparison between preconditioners for PCG. We report precomputing time, total time (including the precompute time) for each precision level, and the corresponding PCG iterations (in parenthesis). The best value in each category is in bold. $\downarrow$: the lower the better.}
\label{tab:comp-x}
\end{table*}

\begin{table*}[htb]
\resizebox{\textwidth}{!}{%
        \scriptsize
\centering
\begin{tabular}{lcccccccc}

\toprule
Task & Method    & Precompute time 	$\downarrow$	& time (iter.) $\downarrow$ 		& time (iter.) 	$\downarrow$	& time (iter.) 	$\downarrow$	& time (iter.) $\downarrow$		& time (iter.) 	$\downarrow$	& time (iter.)$\downarrow$ \\ &

& (s) 		&  until 1e-2 		&  until 1e-4 		&  until 1e-6 		&  until 1e-8 		&  until 1e-10 		&  until 1e-12 \\

\midrule
\multirow{2}{*}{wave-2d}

& AMG & 0.0753 &  0.44  (15) &  0.527  (19) &  0.615  (23) &  0.7  (27) &  0.781  (30) &  0.869  (34) \\

&Ours & 0.0147 & \textbf{0.081} (0) & \textbf{0.081} (0) & \textbf{0.081} (0) & \textbf{0.100} (3) & \textbf{0.156} (12) & \textbf{0.211} (21) \\

\midrule

\multirow{2}{*}{possion-3d}

& AMG & 4.1097 &  5.877  (0) &  10.045  (4) &  13.733  (8) &  17.399  (11) &  20.954  (14) &  24.479  (18) \\
&Ours & 0.4137 & \textbf{3.010} (0) & \textbf{3.220} (2) & \textbf{4.815} (13) & \textbf{6.908} (28) & \textbf{8.749} (41) & \textbf{10.406} (53) \\

\bottomrule

\end{tabular}
}
\caption{Comparison between the Algebraic MultiGrid preconditioner (AMG) and Our proposed preconditioner. We report precompute time, total time (including the precompute time) for each precision level, and PCG iterations (in parenthesis). $\downarrow$: the lower the better.}
\label{tab:amg}
\end{table*}

\subsection{Condition Number Comparison}
\label{supp:cond}

We show the condition number comparison in Table~\ref{supp:table:cond-number}. We observe that IC(2) is the most powerful approach in reducing the condition number of the system matrix $A$, but it is the most computationally expensive approach to derive. Jacobi Method is the least powerful approach in reducing the condition number, and thus the CG iteration, as shown in the Table. This clearly presents the speed-accuracy trade-off. Our method is relatively fast to compute, as shown in~Table.~\ref{tab:pcg}, and also reduces condition number by a relatively large amount. The results reflect and explain that we outperform other numerical baseline methods in total time.  

\begin{table}[htb]
    \centering
    \begin{scriptsize}
    \begin{tabular}{c|cccc}
    \toprule
        Method & Wave-2d & Poisson-2d & Heat-2d & Poisson-3d  \\
        \midrule
        A (original system) & 540272.25 & 43658.16 & 181.56 & 1008.79 \\
        Jacobi & 96.35 & 18712.16 & 165.71 & 225.86 \\
        Gauss-Seidel& 25.20 & 13902.30 & 117.94 & 168.75\\
        IC &  23.22 & 5662.48 & 42.83 & 43.23 \\
        IC(2) & 21.37 & 3742.01 & 34.05  & 37.04\\
        Ours& 23.89 & 8384.31 & 64.23 & 136.86 \\
        \bottomrule
    \end{tabular}
    \caption{{Condition number comparison between various methods}}
    \label{supp:table:cond-number}
    \end{scriptsize}
\end{table}

\subsection{Discussion and Comparison with MultiGrid Preconditioners} 
\label{supp:amg}
We compare it with the algebraic multigrid method (AMG). We adopt the implementation of the commonly used open-source package AMGCL~\cite{Demidov2019}. Results are shown in~\ref{tab:amg}.  We can see from the table that AMG does not perform well on non-elliptic PDEs, such as the hyperbolic PDE (wave-2d). AMG results in the largest number of CG iterations as compared with other baseline numerical approaches, Jacobi, Gauss-Seidel, IC, and IC(2), demonstrating that the multigrid approach is not designed to be general purpose. We also observe that AMG improves the CG iteration number for the elliptic PDE and results in the best CG iteration, as shown in the comparison on Poisson-3d. However, AMG is more computationally expensive and less parallelizable compared to our proposed approach, and thus the derivation time is 10 times longer than our approach. Therefore, AMG takes more total time compared to our proposed approach.

\subsection{Generalizability Comparisons with Learning Physics Simulation Works}\label{supp:exp:network}

\begin{figure*}
    \centering
    \includegraphics[width=\textwidth]{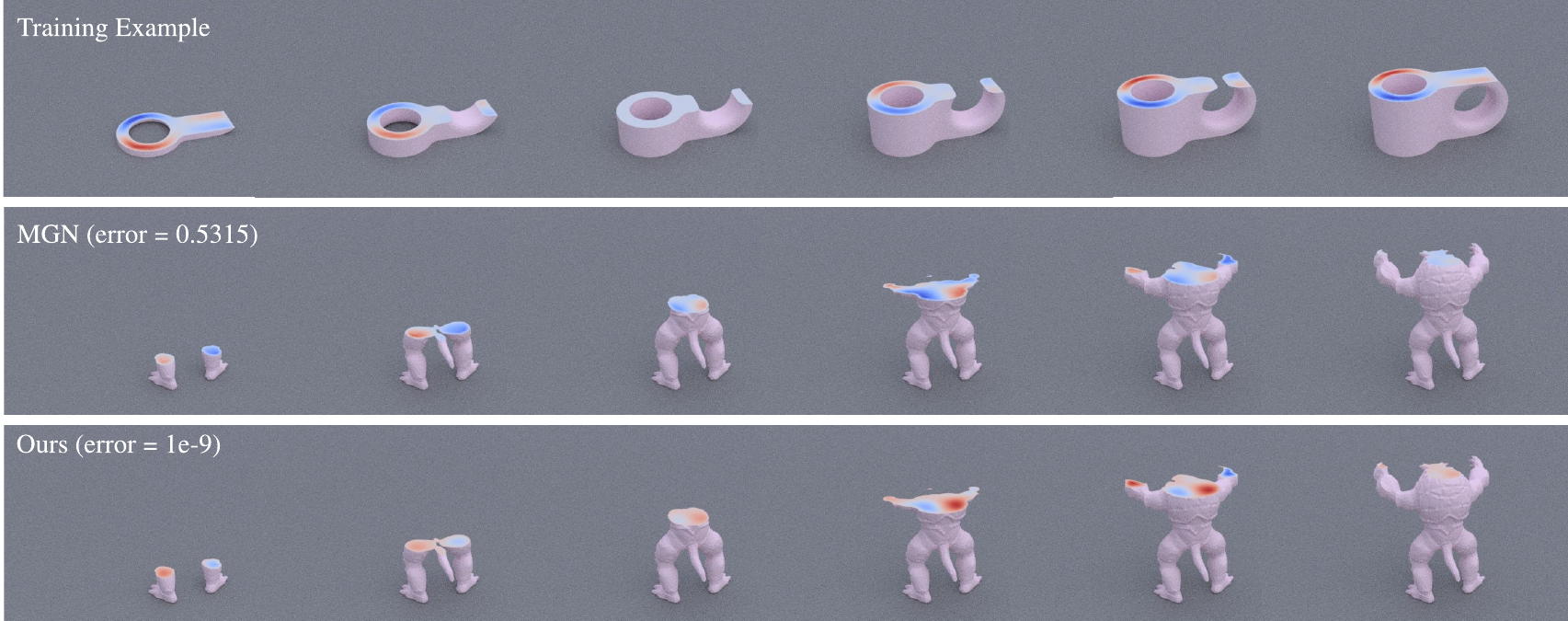}
    \caption{
    {Poisson-3d} on unseen mesh: MGN (middle), our method (bottom), training mesh (top). Left to right: solution fields at different cross-section heights.}
    \vspace{-1mm}
    \label{fig:mesh}
\end{figure*}

We also compare the generalizability of our proposed approach with the learning physics simulation work, MeshGraphNet (MGN) ~\cite{pfaff2020learning}. Both our method and MeshGraphNet are trained on connector-shaped mesh and tested on armadillo mesh, as shown in Fig.~\ref{fig:env}. The results are shown in Fig. ~\ref{fig:mesh}. By comparing the bottom and middle rows in Fig.~\ref{fig:mesh}, we first see that the end-to-end network method (MGN) struggles to generate accurate solutions when deployed on the unseen mesh ($0.5315$ error), whereas our approach achieves arbitrary accuracy by construction ($1e-9$ error threshold here). {This result reflects leveraging neural networks to obtain fast but inaccurate solutions can have more benefit when the solution is not taken as given but used in the context of traditional approaches can have more benefit. 
}

\subsection{Error Accunumation Comparisons with Learning Physics Simulation Works}\label{supp:network-accu}

\begin{figure}[htb]
    \centering
    \includegraphics[width=\linewidth]{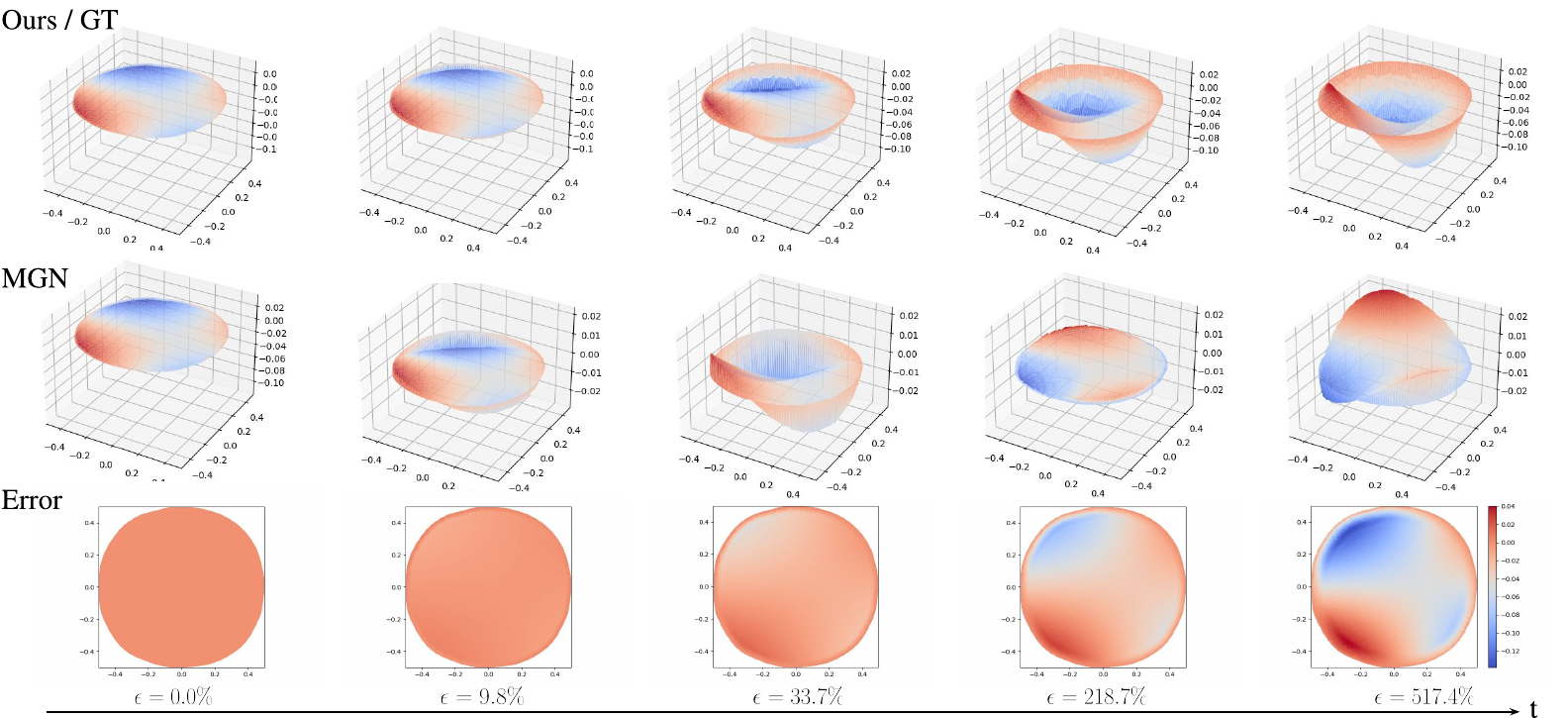}
    \caption{MGN Error accumulation. Field values vs. time step ({wave-2d}): the ground-truth field solved using PCG with 1e-10 convergence threshold (top), the field predicted by MGN (middle), and their difference (bottom), all evaluated at time step 1, 5, 30, 70, 100 (left to right).} 
    \label{fig:acc_error}
    \vspace{-6mm}
\end{figure}

In this experiment, we demonstrate the advantage of our approach over end-to-end network methods: We ensure accurate solutions while end-to-end networks accumulate errors over time. To show this quantitatively, we solve the wave-2d equation for 100 consecutive time steps using our approach and MGN. Fig.~\ref{fig:acc_error} shows that while our approach agrees exactly with the ground truth ($1e-10$ precision), MGN deviates from the ground truth over time. At time step 100, MGN has an error of $517.4\%$. We can expect MGN to be faster in wall-clock time, as MGN only needs to run network inference once, while we need to run network inference to compute the preconditioner followed by running PCG solvers. 

To summarize, MGN is good at estimating solutions rapidly while our approach has the flexibility of achieving arbitrary solution precision, just like a standard CG solver.
Therefore, network methods are suitable for applications where speed dominates accuracy,
while our approach is better for applications that require high precision, e.g., in scientific computing and engineering design.

\end{document}